\documentclass[final]{siamltex}         

\usepackage{graphicx}

\usepackage{subfigure} 

\graphicspath{{C:/NS/Pictures/}}

\title{Research of convergence of the iterative method for solution of the Cauchy problem for the Navier - Stokes equations based on estimated formula.}                                       
\author{ A. Tsionskiy, M. Tsionskiy  }               

\newcommand{\sz}[1]{\mbox{\Large #1 \normalsize}}

\begin{document}             

\maketitle                   

\begin{abstract}
Solution of the Navier-Stokes equations with initial conditions (Cauchy problem) for 2D and 3D cases was obtained in the convergence series form by iterative method using Fourier and Laplace transforms in paper $\cite{TT02}$. For several combinations of problem parameters numerical results were obtained and presented as graphs. 
\nonumber\\

Estimated formula for the border of the parameter area of convergence of the iterative method was obtained in paper $\cite{TT03}$.
\nonumber\\

This paper describes numerical proof of convergence of the iterative method for solution of the Cauchy problem for the Navier - Stokes equations. Usage of estimated formula for the border of the parameter area of convergence of the iterative method is shown for wide ranges of the problem's parameters.

\end{abstract}



\pagestyle{myheadings}
\thispagestyle{plain}
\markboth{A. TSIONSKIY, M. TSIONSKIY}{RESEARCH OF CONVERGENCE OF THE ITERATIVE METHOD}

\section{Introduction}\ 

Example of the solution of the Navier-Stokes equations with initial conditions (Cauchy problem) with a particular radial applied force for 2D case was obtained in paper $\cite{TT02}$. 

We have initial conditions

\begin{equation}\label{eqn400}
\vec{u}(x,0)\; = \; \vec{u}^{0}(x)\; = \;0\;\;\;\;\;\;\;\; (x\in R^{2})
\end{equation}

and a radial applied force

\begin{equation}\label{eqn228}
f_{1\tilde r}(\tilde r, t, \tau, \varphi) = F_{n}\tilde r^{n+1}\sz{e}^{-\mu_{n}^2\tilde r^2}\frac{\sz{e}^{in\tilde \varphi}}{[4\mu_{n}^2\nu (t-\tau)+ 1]^2}
\end{equation}

\begin{equation}\label{eqn228a}
f_{1\tilde \varphi}(\tilde r,\tilde \varphi,\tau) \equiv 0
\quad\quad\quad\quad\quad\quad\quad\quad\quad\quad\quad
\end{equation}

Here ($\tilde r, \varphi$) are the polar coordinates; n is a circumferential mode, n=1,2,3,...; $F_{n}, \mu_{n}$ - constants; \(\nu\) is a positive coefficient of the viscosity; t, \(\tau\)  - time.

A proof of convergence of the iterative method for solution of the Cauchy problem for the 3D Navier - Stokes equations was described in paper $\cite{TT03}$. As one of the results of this proof an estimated formula for the border of the parameter area of convergence of the iterative method was also obtained:

\begin{equation}\label{eqn170x}
\frac{F}{\mu^{4}\nu} < 1
\end{equation}

It is easy to show for an estimate solution of the Cauchy problem for the 2D Navier - Stokes equations that estimated formulas for the border of the parameter area of convergence of the iterative method for different circumferential modes have the same form: 

\begin{equation}\label{eqn170v}
\frac{F_n}{\mu^{4}\nu} < 1\;\;\;\;\;\;\;n=1,2,3,...
\end{equation}

\section{Numerical analysis based on estimated formula}\ 

Numerical results for circumferential modes n = 1, 2, 3, 4, 5 with a radial applied force (\ref{eqn228}) are shown on FIG. 1 - 5. 

We plot curves of the form:

\begin{equation}\label{eqn170w}
\frac{F_n}{\mu^{4}\nu} = 1
\end{equation}

in coordinates $F_n$ (horizontal, logarithmic scale) and $\mu$ (vertical) for various $\nu$. 

Plots for $\nu = 0.01,\; 0.1,\; 0.3,\; 0.75,\; 1.0,\; 1.5\;$ are displayed on FIG. 1.a, 2.a, 3.a, 4.a, 5.a.

Areas of convergence of iterative method are above these  curves.

For each circumferential mode n results are calculated for amplitudes of $F_n$ taken by formula

\begin{equation}\label{eqn170p}
F_n = 10^{k}/n\;\;\;\;\;k = 0,1,2,3.
\end{equation}

Amplitudes of $F_n$ are selected from range $0.2 \leq F_n \leq 1000$.

Since areas of convergence of the iterative method are above displayed curves, we have selected the highest curve from our set of tests for further analysis. 
For all circumferential modes n and all chosen values of $F_n$ the corresponding values of $\mu$ are appeared to be in a range $2.0 \leq \mu \leq 18.0$. All of them are shown as black dots on the curve for $\nu$ = 0.01.
\nonumber\\

For selected values of $F_n$, $\mu$ and $\nu$ comparison of calculated Maximums of Velocity on the first step $\mid\vec u_{1}\mid$ (red plots) and the second step $\mid\vec u_{2}^{*}\mid$ (blue plots) of iterative process is displayed on FIG 1.b, 2.b, 3.b, 4.b, 5.b. The upper graph is shown in coordinates $F_n$ (horizontal, logarithmic scale) and Amplitude of Velocity (vertical, logarithmic scale) for $\nu$ = 0.01. The bottom four graphs are plotted in coordinates $\nu$ (horizontal) and Amplitude of Velocity (vertical, logarithmic scale). For all graphs we see that $\mid\vec{u_{2}^{*}}\mid\; <<\; \mid\vec{u_{1}}\mid$. On FIG 1.c, 2.c, 3.c, 4.c, 5.c the comparison of first step amplitudes and second step amplitudes is displayed in plane $\varphi$ = [0, $\pi$] for various radii r in a range 0 $\leq$ r $\leq$ 3. Also from these graphs we see that for all circumferential modes amplitudes of velocity on the second step of iterative process are smaller than corresponding amplitudes on the first step. Left graphs have full amplitude scale and right graphs are zoomed into low
amplitude values for better display of second step amplitudes. More detailed results of amplitude comparison are shown in $\cite{TT02}$.
\nonumber\\

From all the graphs below we may conclude that the iterative process described in $\cite{TT02}$ and $\cite{TT03}$ has very good convergence. For higher values of $\mu$ with the same values of $F_n$ the convergence of the method will be even better.
\nonumber\\
\nonumber\\
\nonumber\\
\nonumber\\
\nonumber\\
\nonumber\\
\nonumber\\
\nonumber\\
\nonumber\\
\nonumber\\
\nonumber\\
\nonumber\\

   \begin{center}   
     \includegraphics[height=55mm]{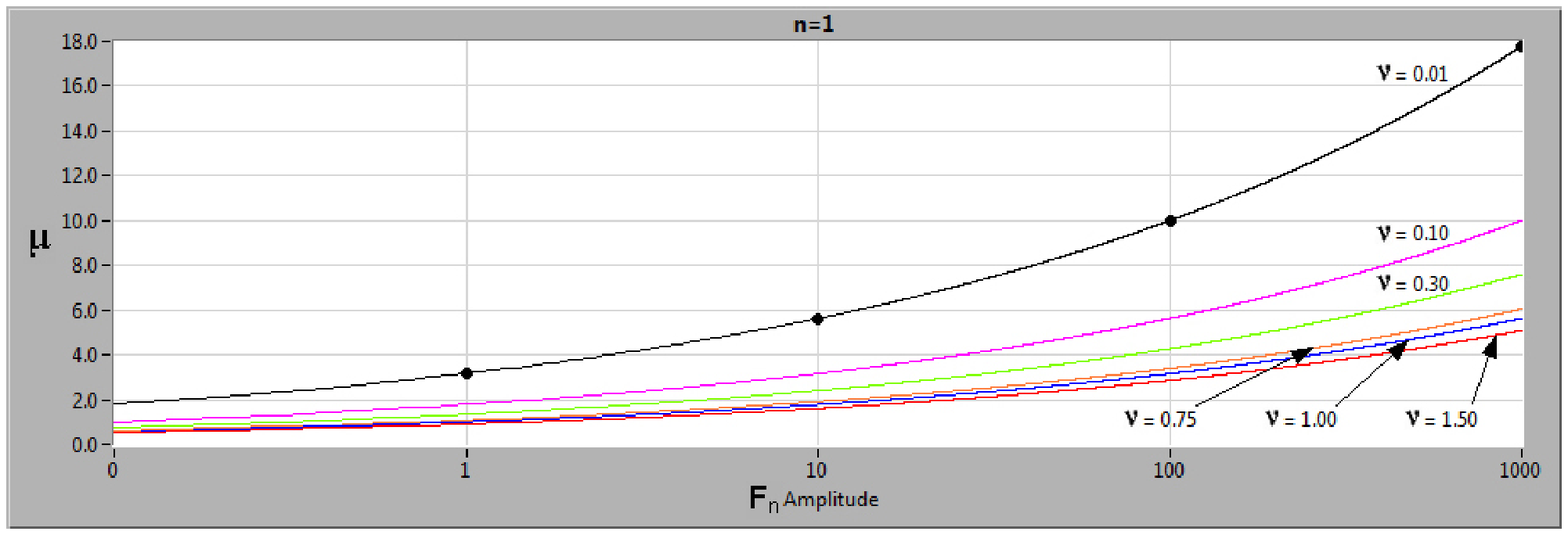}\\
     FIG.1.a.
   \end{center}
   \begin{center}   
     \includegraphics[height=70mm]{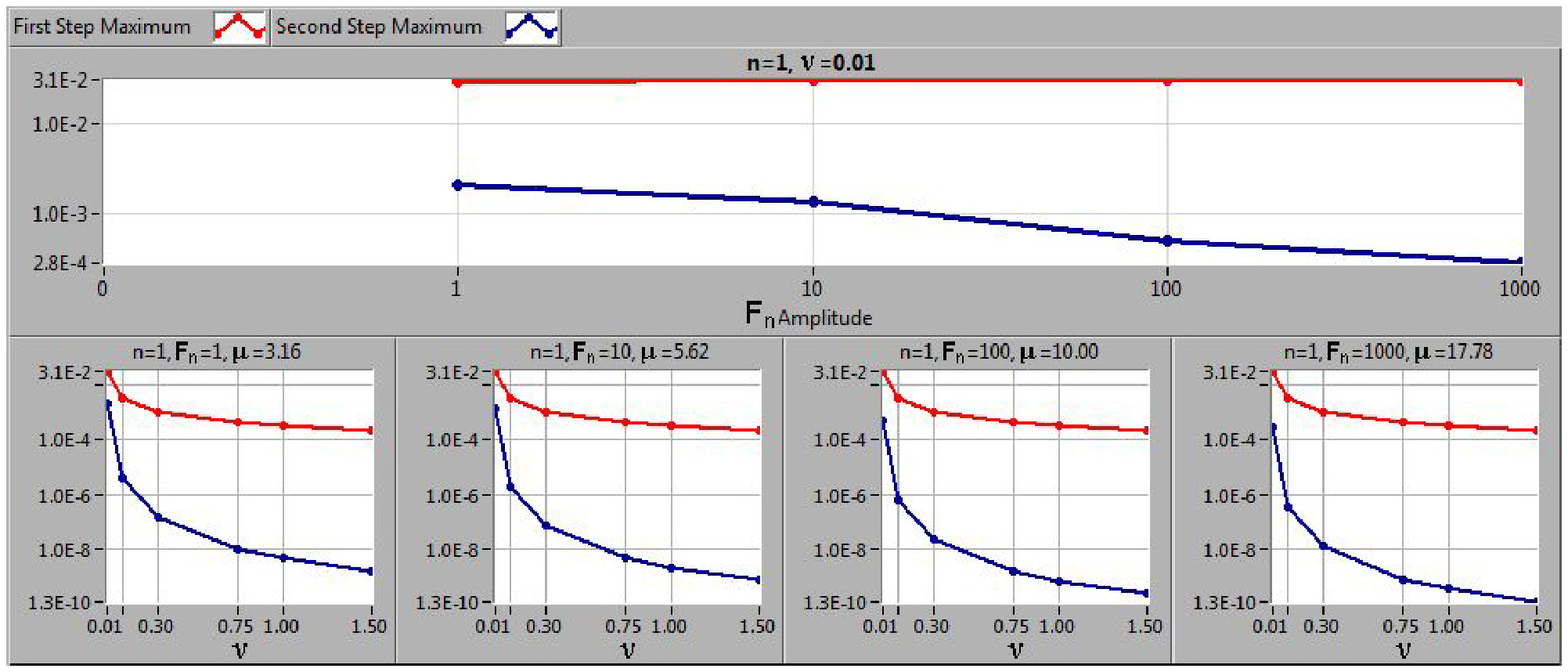}\\
     FIG.1.b.
   \end{center}
   \begin{center}   
     \includegraphics[height=55mm]{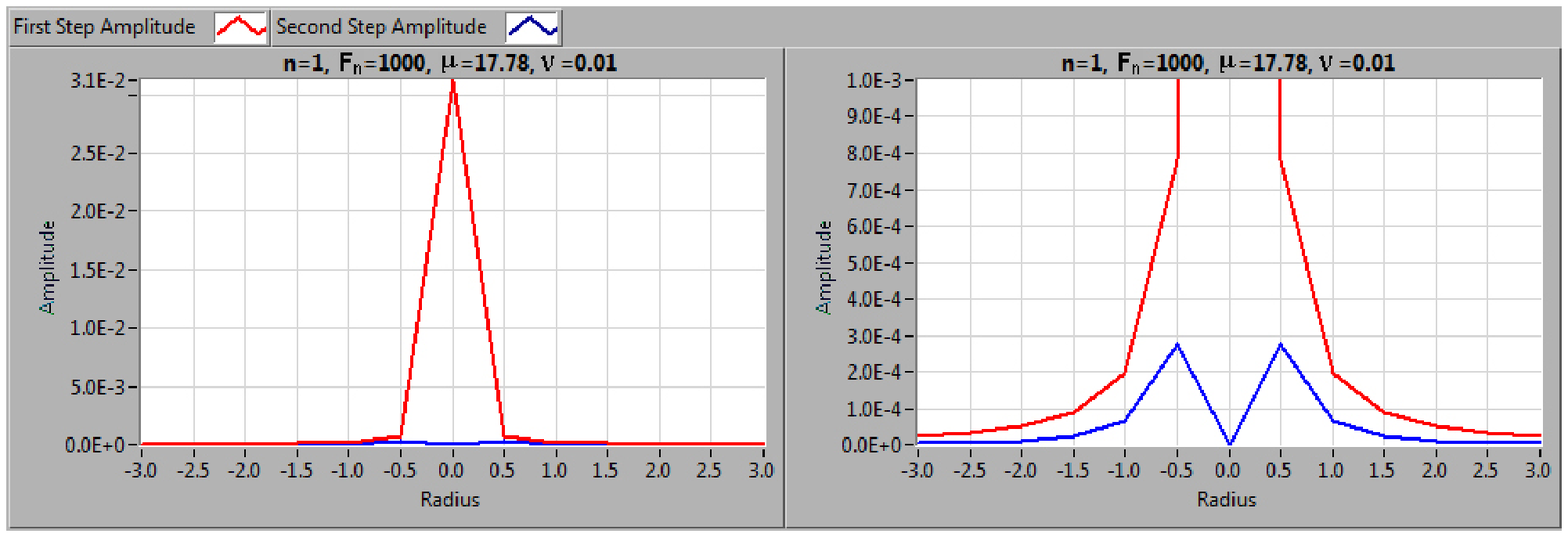}\\
     FIG.1.c.
   \end{center}
   \begin{center}   
     \includegraphics[height=55mm]{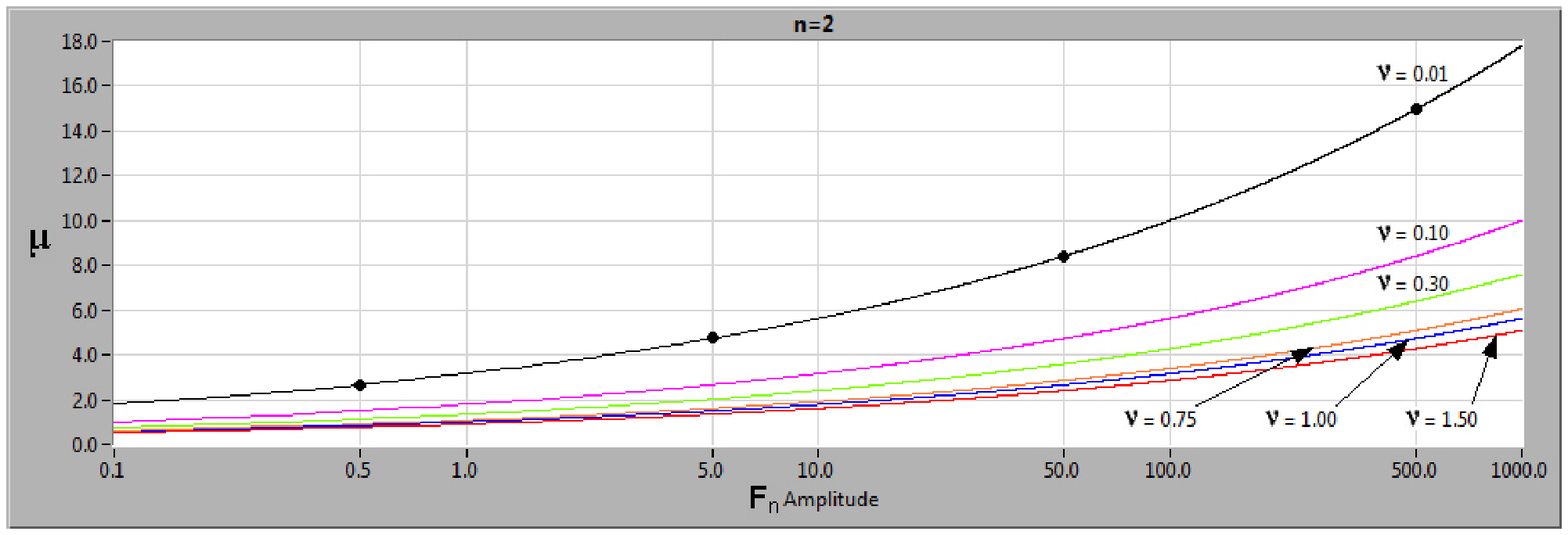}\\
     FIG.2.a.
   \end{center}
   \begin{center}   
     \includegraphics[height=70mm]{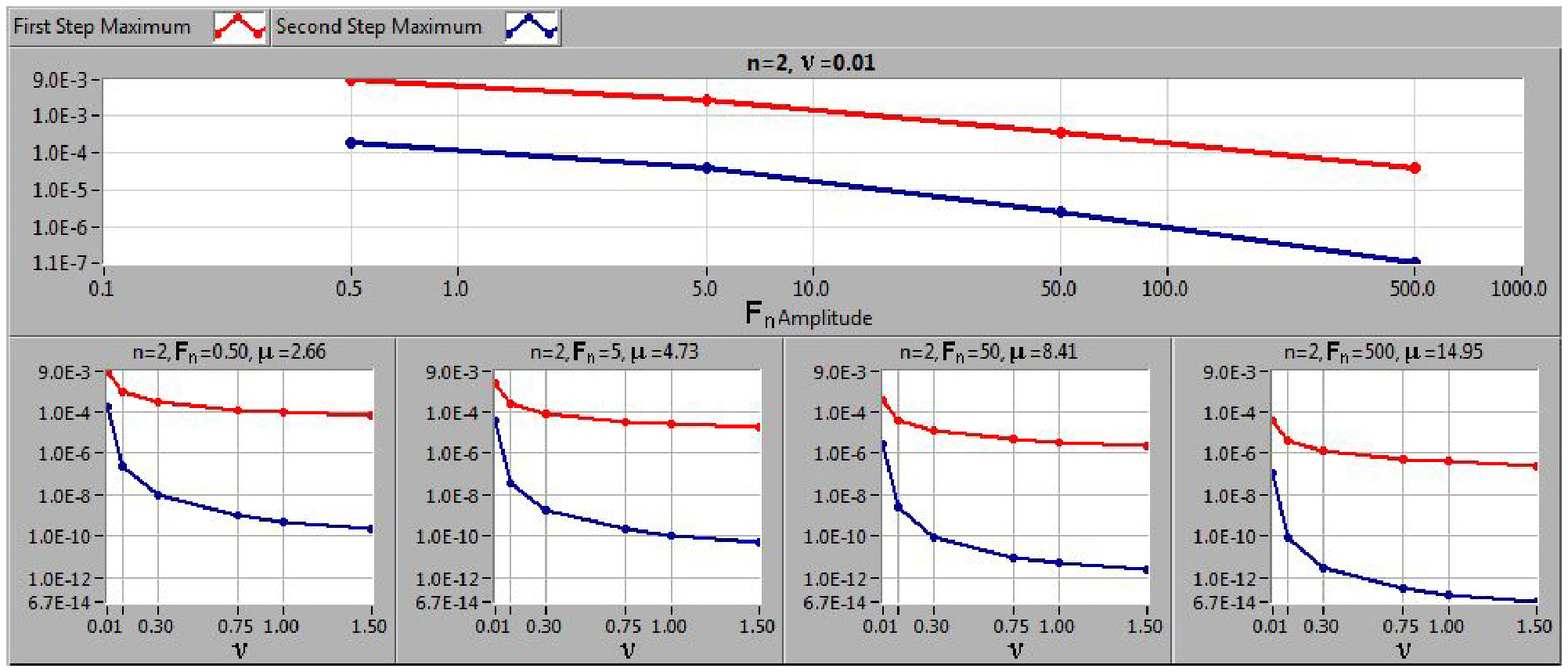}\\
     FIG.2.b.
   \end{center}
   \begin{center}   
     \includegraphics[height=55mm]{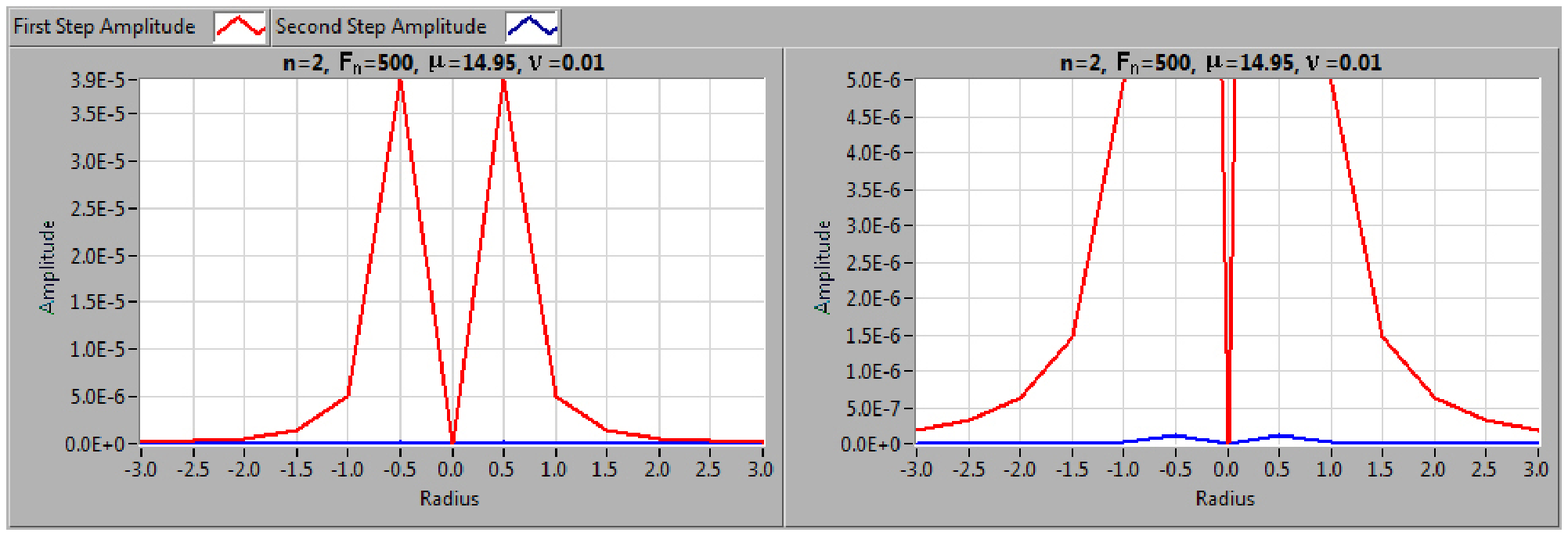}\\
     FIG.2.c.
   \end{center}
   \begin{center}   
     \includegraphics[height=55mm]{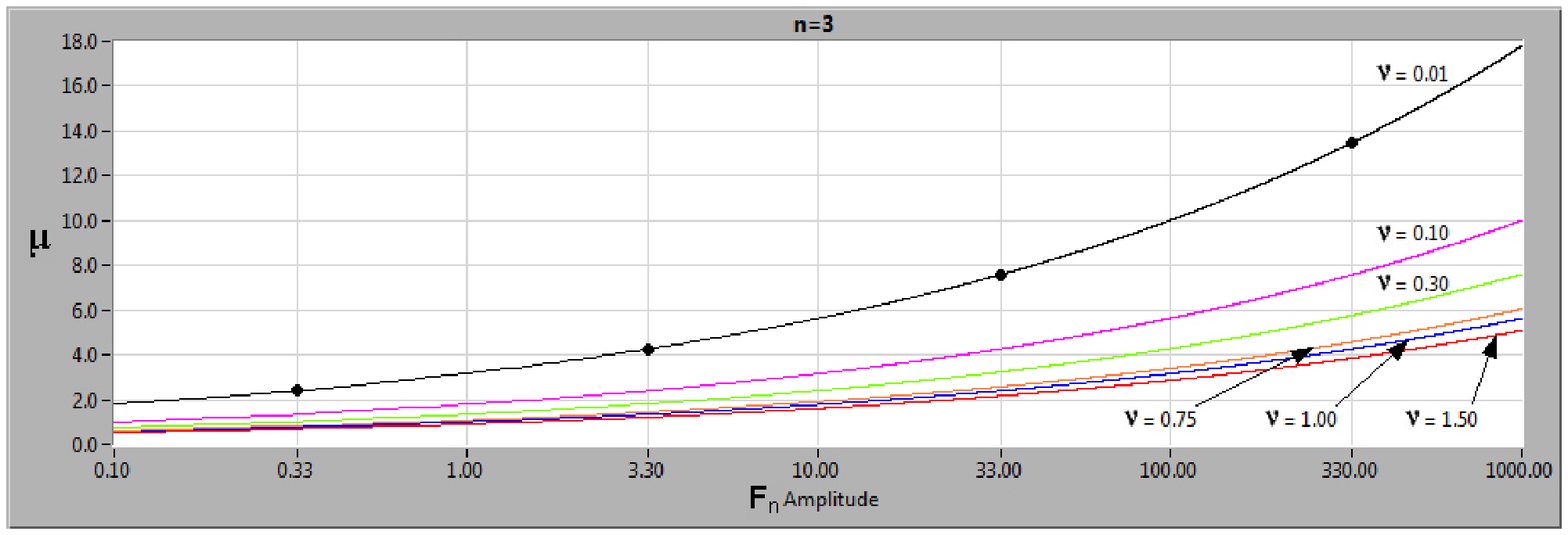}\\
     FIG.3.a.
   \end{center}
   \begin{center}   
     \includegraphics[height=70mm]{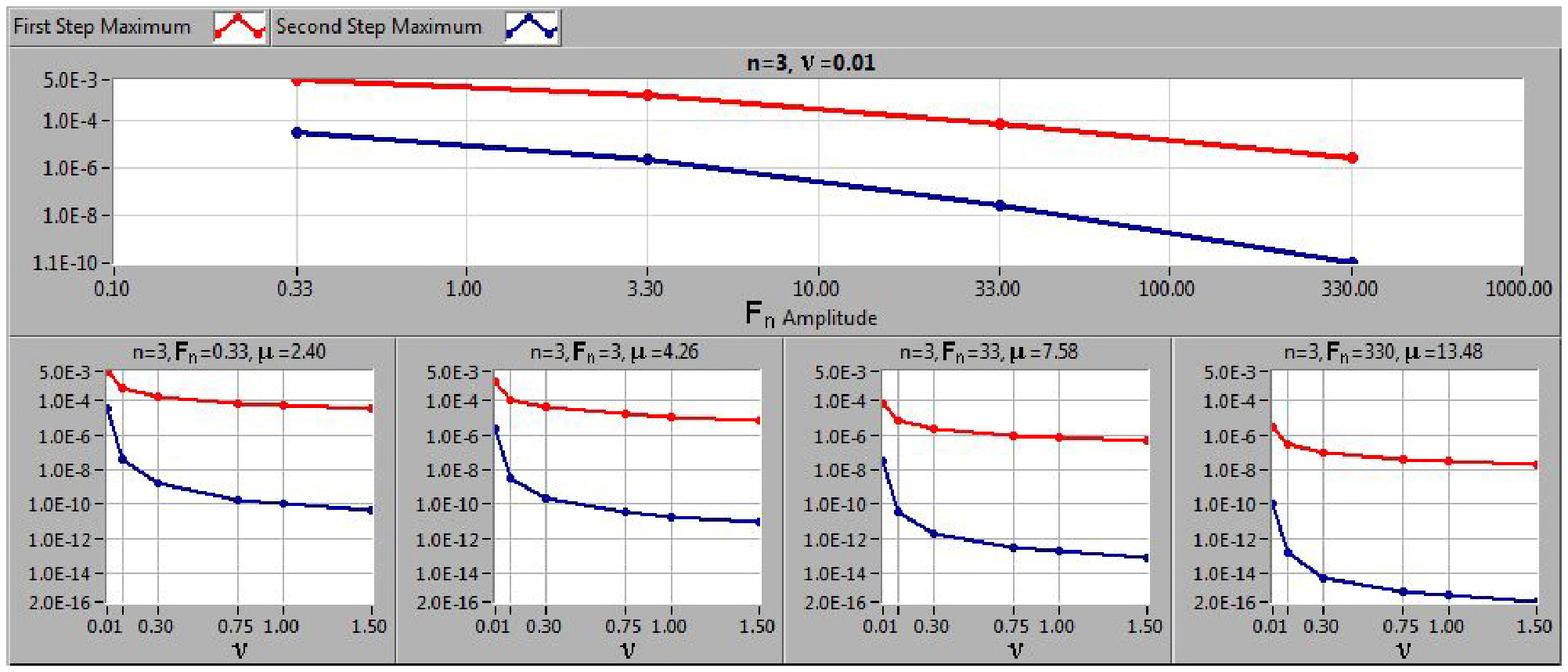}\\
     FIG.3.b.
   \end{center}
   \begin{center}   
     \includegraphics[height=55mm]{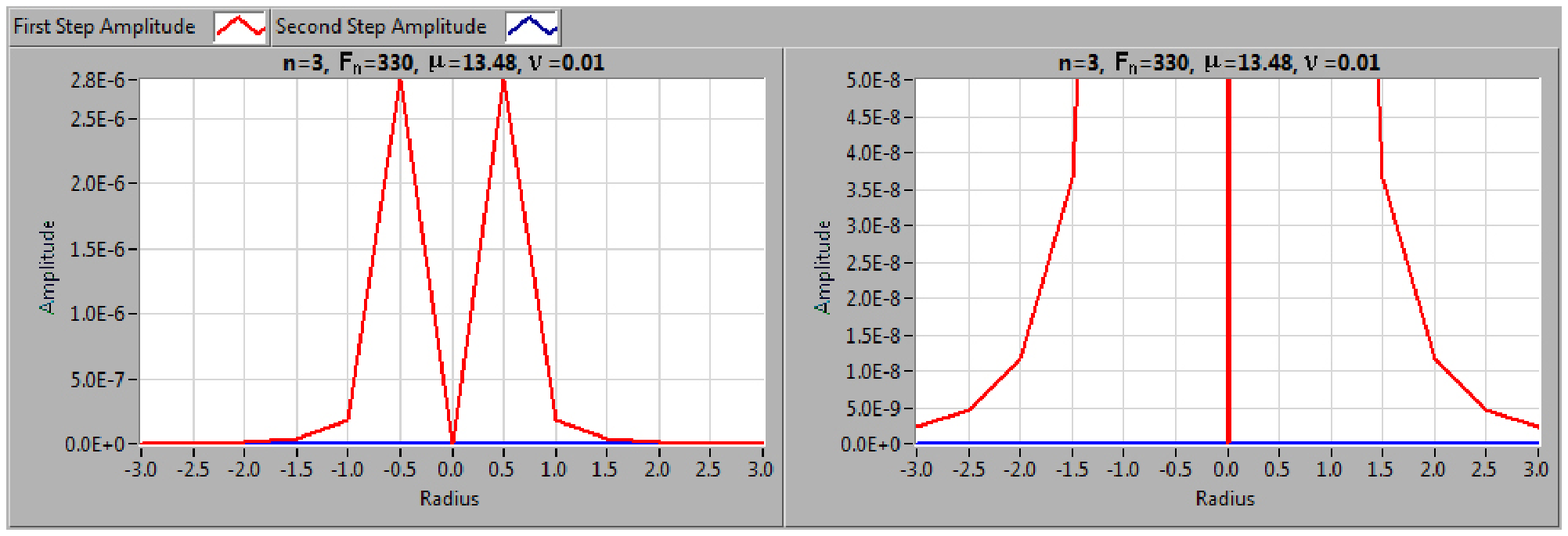}\\
     FIG.3.c.
   \end{center}
   \begin{center}   
     \includegraphics[height=55mm]{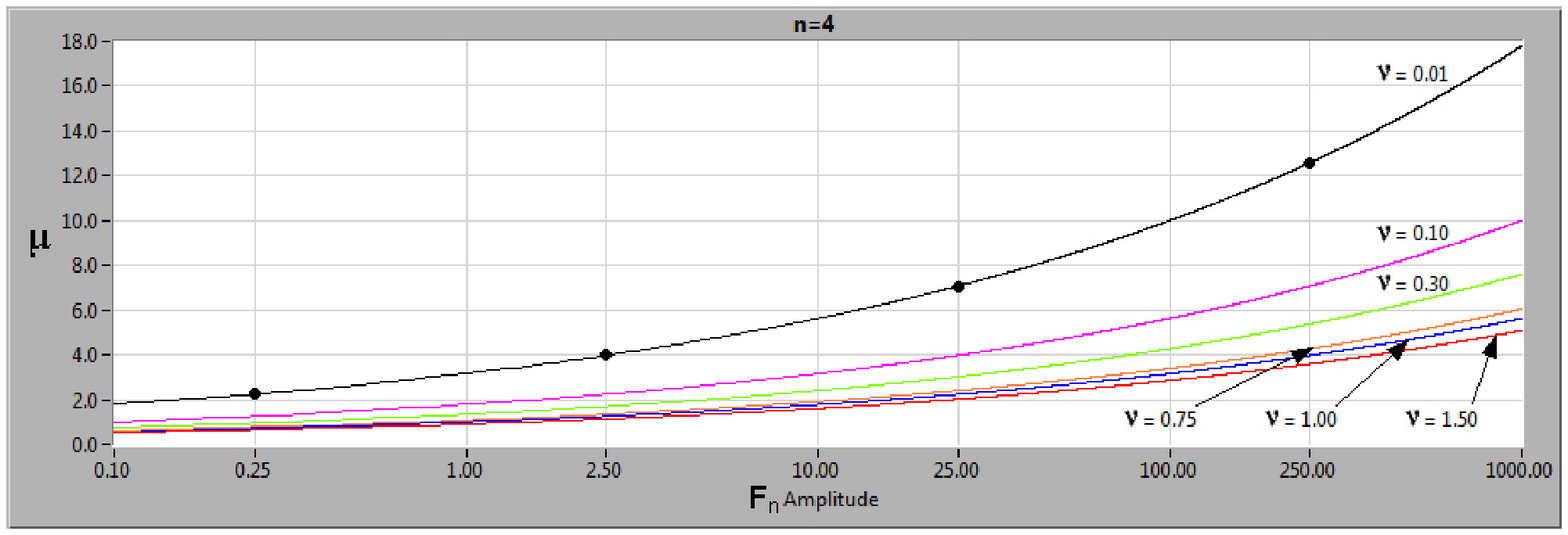}\\
     FIG.4.a.
   \end{center}
   \begin{center}   
     \includegraphics[height=70mm]{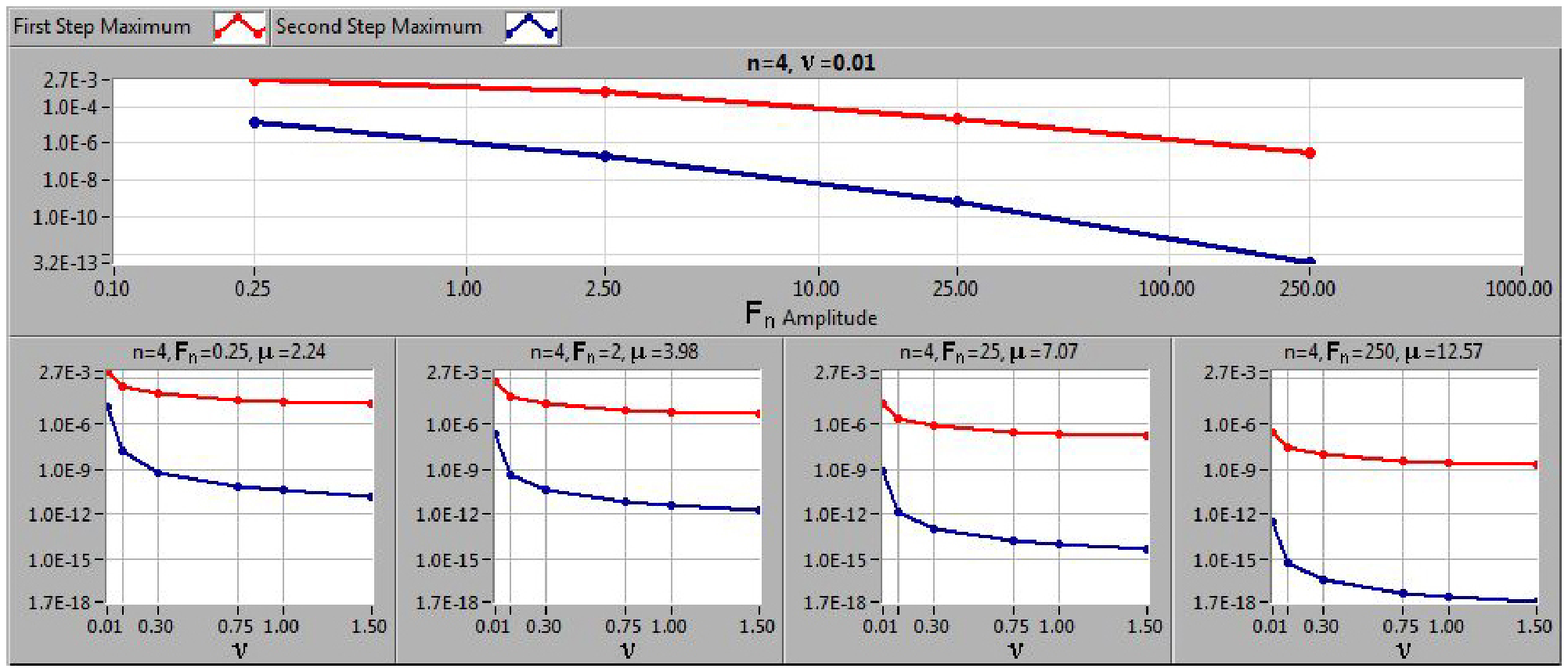}\\
     FIG.4.b.
   \end{center}
   \begin{center}   
     \includegraphics[height=55mm]{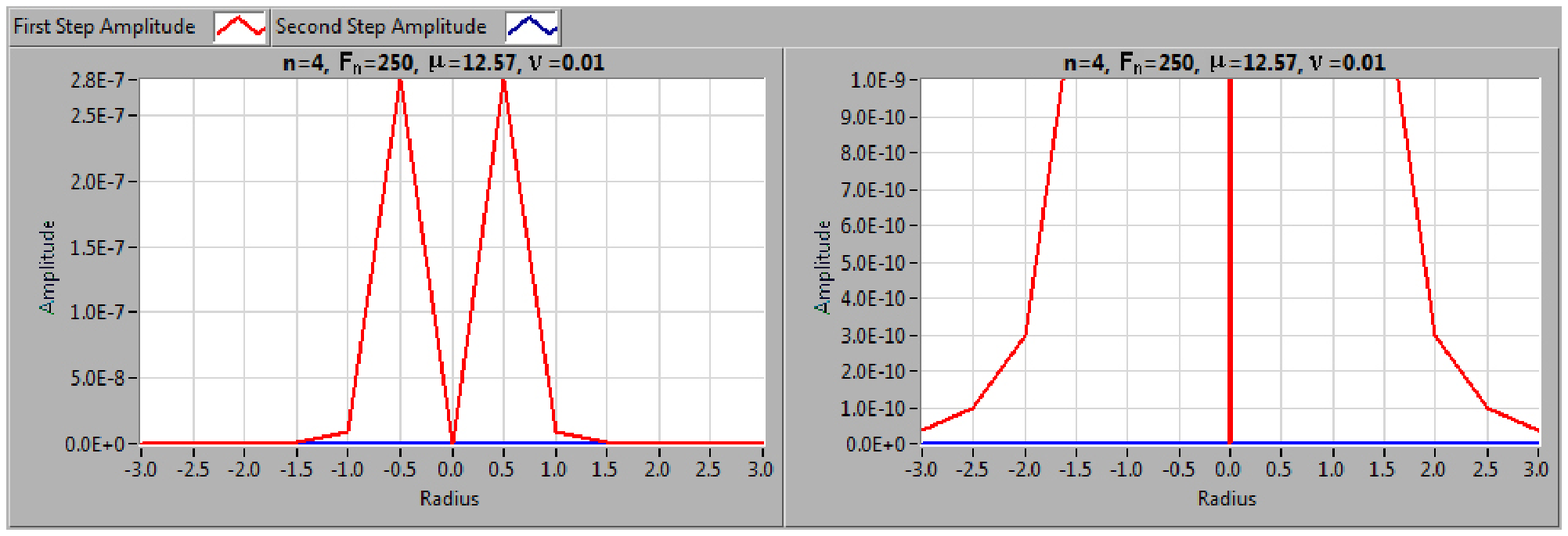}\\
     FIG.4.c.
   \end{center}
   \begin{center}   
     \includegraphics[height=55mm]{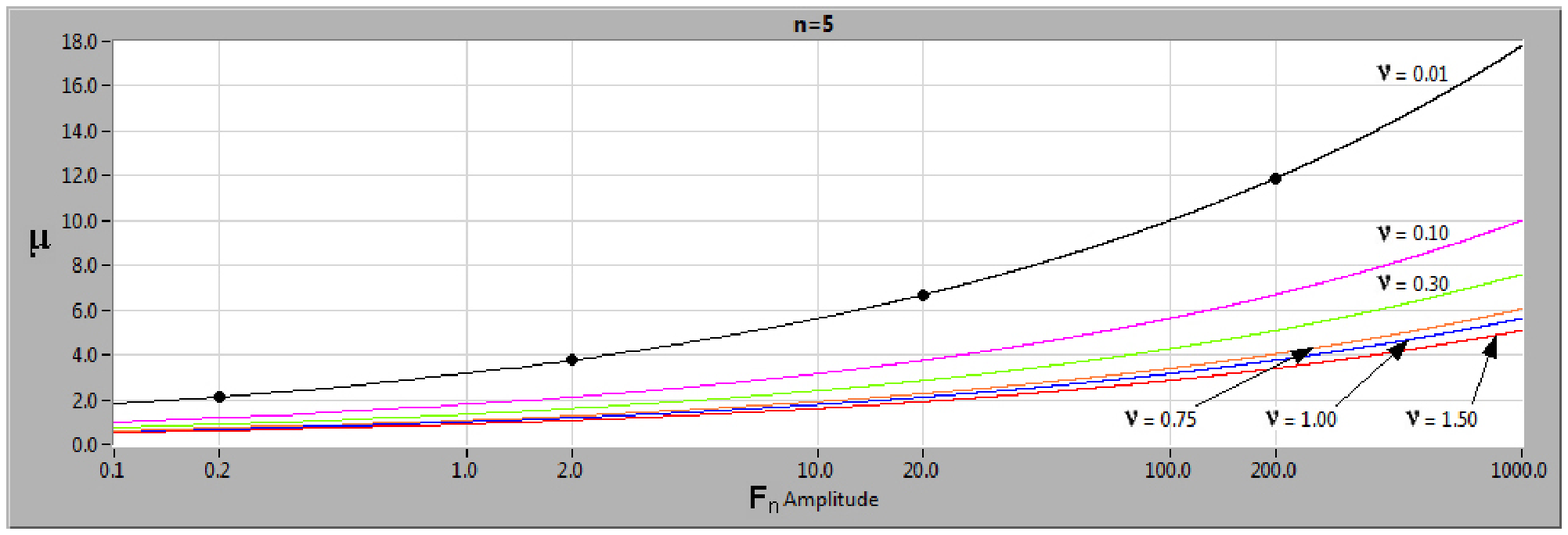}\\
     FIG.5.a.
   \end{center}
$\;$

   \begin{center}   
     \includegraphics[height=70mm]{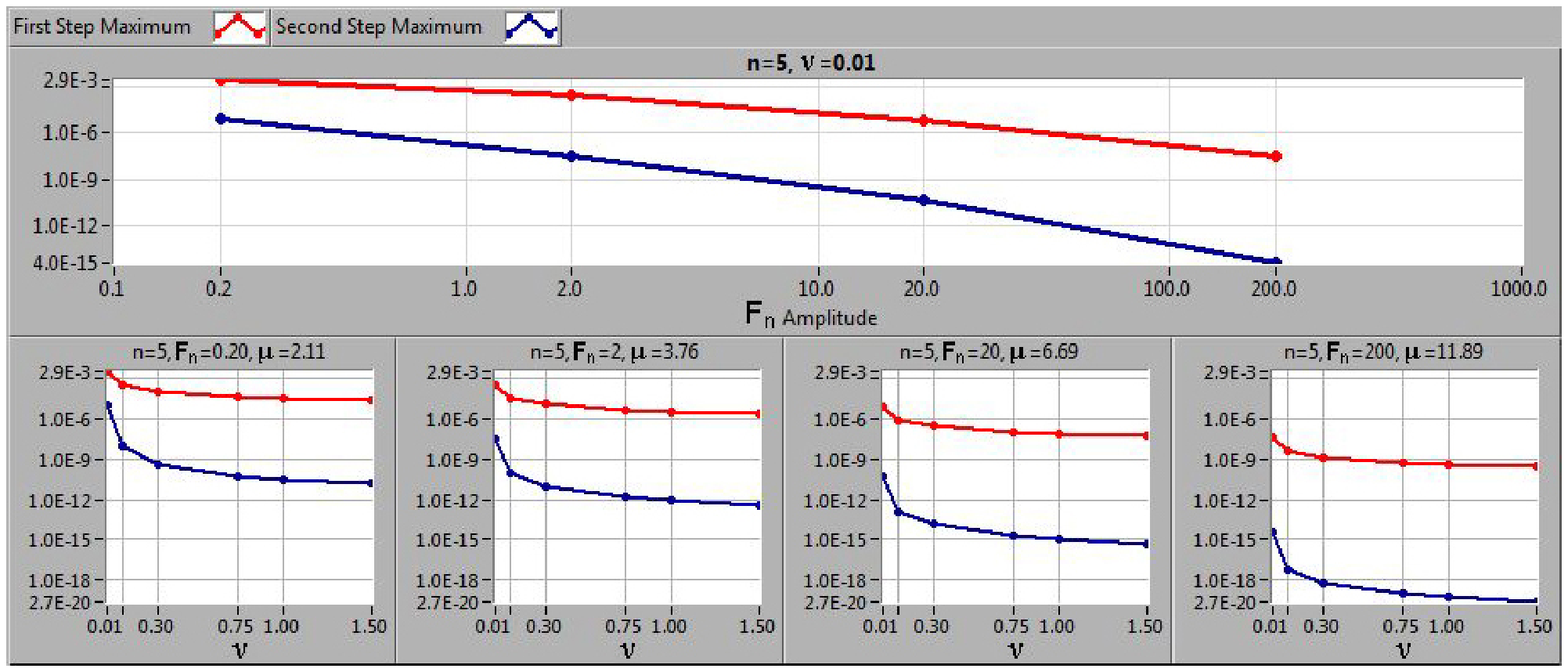}\\
     FIG.5.b.
   \end{center}
$\;$

   \begin{center}   
     \includegraphics[height=55mm]{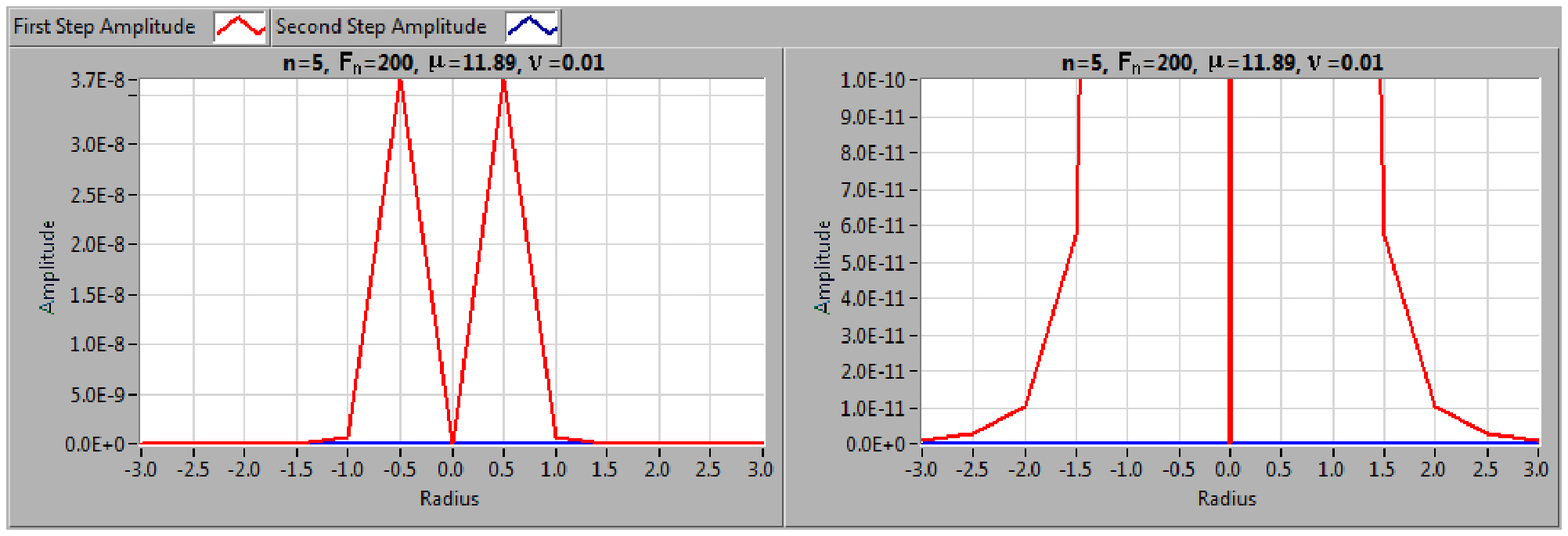}\\
     FIG.5.c.
   \end{center}

\end{document}